\documentclass{article}

\usepackage{arxiv}

\usepackage[utf8]{inputenc} 
\usepackage[T1]{fontenc}    
\usepackage{url}            
\usepackage{booktabs}       
\usepackage{amsfonts}       
\usepackage{nicefrac}       
\usepackage{microtype}      
\usepackage{lipsum}		
\usepackage{graphicx}
\usepackage{natbib}
\usepackage{doi}

\usepackage{graphicx}
\usepackage{amsmath,mathrsfs}
\usepackage{indentfirst,amsfonts,amssymb,newlfont}  
\usepackage{booktabs}
\usepackage{float}
\usepackage{graphics}
\usepackage{subfigure}
\usepackage[usenames]{color}

\begin{document}

\title{A new modified highly accurate Laplace-Fourier method for linear neutral delay differential equations} 


\author{Gilbert Kerr \\
	Department of Mathematics, New Mexico Tech, Socorro, New Mexico, USA.
	\And
	Gilberto Gonz\'alez-Parra  \\
Department of Mathematics, New Mexico Tech, Socorro, New Mexico, USA.\\
}

%
%

\date{Received: date / Accepted: date}

\maketitle

\begin{abstract}
In this article, a new modified Laplace-Fourier method is developed in order to obtain the solutions of linear neutral delay differential equations. The proposed method provides a more accurate solution than the one provided by the pure Laplace method and the original Laplace-Fourier method. We develop and show the crucial modifications of the Laplace-Fourier method. As with the original Laplace-Fourier method, the new modified method combines the Laplace transform method with Fourier series theory. All of the beneficial features from the original Laplace-Fourier method are retained. The modified solution still includes a component that accounts for the terms in the tail of the infinite series, allowing one to obtain more accurate solutions. The Laplace-Fourier method requires us to approximate the formula for the residues with an asymptotic expansion. This is essential to enable us to use the Fourier series results that enable us to account for the tail. The improvement is achieved by deriving a new asymptotic expansion which minimizes the error between the actual residues and those which are obtained from this asymptotic expansion. With both the pure Laplace and improved Laplace-Fourier methods increasing the number of terms in the truncated series obviously increases the accuracy. However, with the pure Laplace, this improvement is small. As we shall show, with the improved Laplace-Fourier method the improvement is significantly larger. We show that the convergence rate of the new modified Laplace-Fourier solution has a remarkable order of convergence $O(N^{-3})$. The validity of the modified technique is corroborated by means of illustrative examples. Comparisons of the solutions of the new modified method with those generated by the pure Laplace method and the original/unmodified Laplace-Fourier approach are presented.
\end{abstract}
\keywords{New modified Laplace-Fourier method \and Neutral delay differential equations \and Laplace transform \and Fourier Series \and Cauchy's residue theorem.
}


\section{Introduction}
Delay differential equations (DDEs) play an imulti-delayportant role in many real world processes \cite{erneux2009applied,halanay2020critical,HIVdelayMurray2000,goubault2018inner,rihan2021delay,shampine2006delay,smith2011introduction}. 
Mathematical models based on different types of DDEs have been used to describe and study the dynamics of many infectious diseases \cite{Hethcote,SEISdelayEXP,SISdelayHethcote,SEIRdelayYan,DelayHepB,DelayHeroin,DelayRabiesOneBackSI,HIVdelayMurray2000,gonzalez2022mathematical}. 

The complexity of DDEs is usually greater than those of ordinary differential equations (ODEs). The delays can cause a stable equilibrium to become unstable and periodic solutions can arise \cite{bachar2019periodic,biccer2021periodic,bicer2018asymptotic,breda2014stability,faria2019existence,he2005periodic,SISdelayHethcote}. Recently, there has been an increasing interest in the development of methods for finding and studying solutions of a variety of DDEs including some of fractional order \cite{biccer2021periodic,castro2019exact,el2019analytical,garcia2018exact,gumgum2020lucas,jamilla2020explicit,jamilla2020solutions,khodabandehlo2022numerical,kumar2022efficient,luo2022semianalytical,mohammadian2022numerical,panghal2022neural,senu2022numerical}. The classical method for solving DDEs is the method of steps (MoS), which generates an analytical piecewise solution. This method is based on a stepwise process that requires solving a differential equation in different time intervals  \cite{bauer2013solving,kalmar2009stability,kaslik2012analytical,smith2011introduction,heffernan2006solving}. However, for many DDEs, the MoS encounters the phenomenon of expression swell \cite{bauer2013solving,heffernan2006solving,kerr2022accuracy}. Thus, for many DDEs the MoS is not able to generate the solution for the whole time domain. However, in some cases the MoS can be used to find a closed-form solution \cite{bauer2013solving,heffernan2006solving,kerr2022accuracy}. 

There are numerous classifications for the DDEs, amongst which one important type are 
the neutral delay differential equations (NDDEs). In these equations the time delays appear in the state derivative \cite{bazighifan2021neutral,fabiano2018spline,fabiano2013semidiscrete,faheem2021collocation,liu2019stability,qin2019continuous,philos2001periodic,raza2019haar,santra2020oscillation}. NDDEs arise in numerous mathematical models  involving problems in engineering and science \cite{gu2003survey,kim2019feedback,kuang2012delay,xu2017balancing}. NDDEs have been solved by a variety of methods  \cite{enright1998convergence,jamilla2020solutions,qin2019continuous,xu2016local}.

The main aim of this paper is to present a new modified Laplace-Fourier method, which produces an essentially analytical solution that significantly reduces the error. The novel original Laplace-Fourier method is based on the combination of the Laplace transform and Fourier series theory to obtain the solutions of linear NDDEs \cite{kerr2022new}. The solution is composed of a closed-form part and an infinite series \cite{kerr2022new}. The original Laplace-Fourier method requires computing the poles of a Laplace transformed equation and computing the residues by Cauchy's residue theorem. An original formula developed by the authors provides an excellent
approximate for the location of the poles. However, in the original paper, the formula which was used for approximating the residues was, in some cases, less accurate. In this paper, the authors develop a new approach that provides a much more accurate approximation for the residues, and therefore more accurate solutions. 


The pure Laplace method has been used recently to solve fractional linear DDEs \cite{golmankhaneh2023initial}. In \cite{kerr2022new} it was shown that a previous Laplace-Fourier method can improve the accuracy of the solution of linear DDEs. The Laplace-Fourier method improves the pure Laplace method by taking into account the contribution of all the terms in the pure Laplace solution (infinite series) for terms such that $N>k$. This is accomplished by approximating the tail of the infinite series with a harmonic Fourier series. Then, in turn, determining an analytic formula for the series. This formula, which constitutes an important component of the solution, is a piecewise continuous function that is composed of polynomials. The Laplace-Fourier method has the versatility to modify the degree of polynomials in order to improve the accuracy. The Laplace-Fourier method also requires us to approximate the formula for the residues with an asymptotic expansion. This is essential to enable us to proceed, and use the Fourier series results for solving linear NDDEs.

%

\section{Definitions for NDDEs and Laplace transform}

To facilitate the introduction of the new modified Laplace-Fourier method we will recap some of the basic
foundations (for more details see \cite{kerr2022accuracy,kerr2022new}). Let us focus on the subclass of NDDE's of the following form
\begin{equation}
y^{\prime }(t)=ay(t)+by^{\prime }(t-\tau )+cy(t-\tau ), t>0, b\neq 0.
\label{NDDE 1}
\end{equation}
For the history function $H(t)$ , $t\in [-\tau ,0]$.
\subsection{Applying of the Laplace transform} 
Taking the Laplace Transform of Eq. (\ref{NDDE 1}) and solving for $Y(s)$, one obtains
\begin{equation}
Y(s)=\dfrac{N(s)}{s-a-(bs+c)\,e^{-s\tau }}, \label{alg3}
\end{equation}
where $N(s)=y(0)-b\,y(-\tau )+(bs+c)\,e^{-s\tau }\int_{-\tau }^{0}y(v)\,e^{-sv}dv$.
Assuming that any singularities resulting from the integral term in $N(s)$
are removable, the relevant poles are determined by the roots of the equation
\begin{equation}
D(s)=s-a-(bs+c)\,e^{-s\tau }=0.  \label{dms}
\end{equation}
The exact number of real roots (at most 2) is determined by the parameters $a,b,c$ and the time delay $\tau$. However, for the piecewise polynomial component of the modified Laplace-Fourier solution only the complex poles play a role.

\subsection{Complex Poles}

Dividing the characteristic equation (\ref{dms}) by $s$ reveals that the approximate location of the poles (for $\left\vert s\right\vert $ large and $\mathfrak{Im}(s)>0$) is given by
\begin{equation}
s_{k}=\frac{\ln (b)+2k\pi i}{\tau }\text{ , }k\in \mathbb{N},  \label{sk even}
\end{equation}
when $b>0$. For $b<0$ we have a different formula, but similar (\cite{kerr2022accuracy}).
These formulas can be used as the initial guesses for determining the actual poles (details in \cite{kerr2022accuracy}). 
Applying the inverse Laplace transform to Eq. (\ref{alg3}) and Cauchy's residue theorem one can obtain the solution $y(t)$ of the NDDE (\ref{NDDE 1}). However, this requires finding the residues at the infinitely many poles. Applying L'Hopital's rule, one gets that the related residues at the poles are given by 
\begin{equation}
c_{k}=\dfrac{N(r_{k})}{1+(br_{k}\tau -b+c\tau )\,e^{-r_{k}\tau }},
\label{ck lap}
\end{equation}
where $r_{k},k\in \mathbb{N}$ are the sequence of complex poles above the imaginary axis. Once the poles have been determined we can use the Cauchy
residue theorem to determine $y(t)$. Assuming that there is one real root ($r_{0}$) then we can write the solution as
\begin{equation}
y(t)=c_{0}\,e^{r_{0}t}+2\sum_{k=1}^{\infty }\mathfrak{Re}\left(
c_{k}e^{r_{k}t}\right).  \label{yt lap}
\end{equation}
In the case of two real roots ($r_{0}$ and $r_1$) one obtains
\begin{equation}
y(t)=c_{0}\,e^{r_{0}t}+c_{1}\,e^{r_{1}t}+2\sum_{k=2}^{\infty }\mathfrak{Re}\left(
c_{k}e^{r_{k}t}\right).  \label{yt lap1}
\end{equation}
Thus, regardless of the number of real roots we can always obtain a solution in terms of an infinite non-harmonic series for the linear DDE (\ref{NDDE 1}) if $H(t)\in C^{1}$. Moreover, the solution is unique (see \cite{hale2013introduction,kuang1993delay}). The solution of the NDDE (\ref{NDDE 1}) at some point would have a discontinuous derivative \cite{kuang1993delay}. The smoothing property is not valid for NDDEs and the solution $y(t)$ is always as smooth as the history function $H(t)$ (see \cite{kuang1993delay}). For instance, if the history function $H(t)$ is discontinuous at some point $t=t_0$ then the analytical solution of the NDDE (\ref{NDDE 1}) would be discontinuous at $t=t_0+\tau$.

\subsection{Laplace-Fourier solution}

The main shortcoming of the pure Laplace transform method, when utilized for solving DDE's, is the slow convergence rate of the series solution in the
vicinity of the points $t=m\tau ,m\in \mathbb{N}$ (\cite{kerr2022accuracy}). 
The Laplace-Fourier is able to improve the convergence rate, without increasing the number of terms in the series \cite{kerr2022new}. For the Laplace-Fourier method we shall assume that $H(t)\in C^{\infty}$ in order to be able to guarantee the existence of the Laplace-Fourier solution. Recall, the approximate location of the poles (for $\left\vert s\right\vert$ large, $b>0$ amd $\mathfrak{Im}(s)>0$) are given by Eq. (\ref{sk even}). Therefore, the tail of the series (for $k\geq N$) can be approximated by
\begin{equation}
y_{a}(t)\simeq 2\sum_{k=N}^{\infty }\mathfrak{Re}\left( c_{k}e^{\left[ \frac{\ln
(b)+2k\pi i}{\tau }\right] t}\right) \simeq 2e^{\frac{\ln (b)}{\tau }%
t}\sum_{k=N}^{\infty }\mathfrak{Re}\left( c_{k}e^{\frac{2k\pi i}{\tau }t}\right).
\label{y app 1}
\end{equation}
Observe, that being able to remove the $\ln(b)$ term from inside of the sum leaves one with a regular, harmonic Fourier series (inside the sum). It can also be shown (we will verify this later, in Section \ref{sect3}) that the asymptotic expansion for the $c_{k}$'s at $k=\infty $ is (typically) of the form
\begin{equation}
c_{k}^{a}\simeq \frac{a_{2}}{(\frac{2k\pi i}{\tau })^{2}}+\frac{a_{3}}{(%
\frac{2k\pi i}{\tau })^{3}}+ \cdots,  \label{cka b pos}
\end{equation}
in which all the coefficients $a_{m}$ are real. Hence, labeling ${\alpha }_{k}=\dfrac{2k\pi }{\tau }$ one gets
\begin{equation}
y_{a}(t)\simeq 2e^{\mathfrak{\ln (b)}{\tau }t}\sum_{k=N}^{\infty }\mathfrak{Re}
\left( \left[ \frac{-a_{2}}{{\alpha }_{k}^{2}}+\frac{ia_{3}}{{\alpha }%
_{k}^{3}}+\cdots\right] e^{{\alpha }_{k}it}\right),  \label{y app 2}
\end{equation}
or using Euler's formula one obtains
\begin{equation}
y_{a}(t)\simeq 2e^{\frac{\ln (b)}{\tau }t}\sum_{k=N}^{\infty }\left[ \frac{-a_{2}}{{\alpha }_{k}^{2}}\cos ({\alpha }_{k}t)-\frac{a_{3}}{{\alpha }
_{k}^{3}}\sin ({\alpha }_{k}t)+\cdots\right].  \label{y app 3}
\end{equation}
Which can then be rewritten as 
\begin{equation}
\begin{array}{lll}
y_{a}(t) & \simeq & 2e^{\frac{\ln (b)}{\tau }t}\sum\limits_{k=1}^{\infty }%
\left[ \frac{-a_{2}}{{\alpha }_{k}^{2}}\cos ({\alpha }_{k}t)-\frac{a_{3}}{{%
\alpha }_{k}^{3}}\sin ({\alpha }_{k}t)+\cdots\right] \\ 
&  &  \\ 
&  & -2e^{\frac{\ln (b)}{\tau }t}\sum\limits_{k=1}^{N-1}\left[ \frac{-a_{2}}{%
{\alpha }_{k}^{2}}\cos ({\alpha }_{k}t)-\frac{a_{3}}{{\alpha }_{k}^{3}}\sin({\alpha }_{k}t)+\cdots\right].
\end{array}
\label{yat_sep}
\end{equation}
Observe that the infinite series in the sum which starts at $k=1$ are, by design, an exact match with the classical type harmonic Fourier Series \cite{brown2012fourier,oberhettinger2014fourier}. Therefore, using the relevant Fourier series formulas, enables us to replace this sum with a closed-form analytic expression, in the form of a piecewise continuous function.

\subsection{Fourier Series results}

The results of the previous section were obtained by assuming that $b>0$, where the approximate location of the poles is given by Eq. (\ref{sk even}). For the sake of clarity let us keep this assumption in this section. For the case where $b<0$ the computations are similar and can be found in \cite{kerr2022new}. Some of the following results can also  be found in \cite{kerr2022new}. However, we include them here for self readability.
  
On the interval $0<x<\tau $, with ${\alpha }_{n}=\dfrac{2n\pi }{\tau }$ we
have that
\begin{equation}
\sum\limits_{n=1}^{\infty }{\frac{{\cos (\alpha }_{n}{x)}}{{\alpha }_{n}^{2}}%
=}\frac{3x^{2}-3x\tau +\tau ^{2}/2}{12}=p_{2}(x).  \label{p2}
\end{equation}
Subsequent results or polynomials of higher degree may be obtained via integration \cite{kerr2022new,oberhettinger2014fourier}. 
Using these polynomials enables us re-write $y_{a}(t)$ as
\begin{align}
y_{a}(t)\!& \simeq 2e^{\frac{\ln (b)}{\tau }t}\left[
-a_{2}p_{2}(t)-a_{3}p_{3}(t)+a_{4}p_{4}(t)\right] \notag\\ 
& -2e^{\frac{\ln (b)}{\tau }t}\sum\limits_{k=1}^{N-1}\!\left[ \frac{-a_{2}}{%
{\alpha }_{k}^{2}}\cos({\alpha }_{k}t)\!-\!\frac{a_{3}}{{\alpha }_{k}^{3}}\sin ({\alpha }_{k}t)\!+\!\!\frac{a_{4}}{{\alpha }_{k}^{4}}\cos ({\alpha }_{k}t)\right]\!,0<t<\tau. \label{y app 4}
\end{align}
Therefore 
\begin{align}
y_{a}(t)\!&\simeq 2e^{\frac{\ln (b)}{\tau }t}P_{e}(t)\notag\\
&-2e^{\frac{\ln (b)}{\tau }%
t}\sum_{k=1}^{N-1}\!\left[ \frac{-a_{2}}{{\alpha }_{k}^{2}}\cos({\alpha }_{k}t)\!-\!\frac{a_{3}}{{\alpha }_{k}^{3}}\!\sin({\alpha }_{k}t)\!+\!\frac{a_{4}}{{\alpha }_{k}^{4}}\!\cos({\alpha }_{k}t)\right]\!\!,0<t<\infty.  \label{ya_tail}
\end{align}
Where $P_{e}(t)$ is the $\tau $-periodic extension of the relevant polynomial on the basic interval: $0<t<\tau$. Then, combining the above formula with Eq. (\ref{yt lap}) one gets 
\begin{equation}
y(t)\simeq c_{0}\,e^{r_{0}t}+2e^{\frac{\ln (b)}{\tau }t}P_{e}(t)+2%
\sum_{k=1}^{N}\mathfrak{Re}\left( c_{k}e^{r_{k}t}-c_{k}^{a}e^{\left[ \frac{\ln
(b)+2k\pi i}{\tau }\right] t}\right).  \label{ya_proc}
\end{equation}
in which $c_{k}^{a}$ is given by Eq. (\ref{cka b pos}). We also note that it is possible to extend Eq. (\ref{ya_proc}) by including polynomials of higher degree. 

\subsection{Further insight into the Laplace-Fourier method}

The pure truncated Laplace solution does not include the contribution of all the terms.The most novel aspect of the Laplace Fourier solution is that it does account for these excluded terms \cite{kerr2022new}. 
The Laplace-Fourier method is very accurate in the case where $ab+c=0$. This is due to the fact that the estimates (\ref{sk even}) provide an exact formula
for the poles, since one gets
\begin{equation}
D(s_{k})=s_{k}\!-a-\left(bs_{k}\!-c\right)\,e^{-(\frac{\ln (b)+2k\pi i}{\tau })\tau }  
=s_{k}\!-a\!-\!\frac{1}{b}(bs_{k}\!+\!c)=-\!\left(\!a\!+\!\frac{c}{b}\right)\!\!.\!
\label{Ds_k}
\end{equation}
Thus Eq. (\ref{ya_proc}) becomes
\begin{equation}
y(t)\simeq c_{0}e^{r_{0}t}+2e^{\frac{\ln (b)}{\tau }t}P_{e}(t)+2e^{\frac{\ln
(b)}{\tau }t}\sum_{k=1}^{N}\mathfrak{Re}\left( [c_{k}-c_{k}^{a}]e^{\frac{2k\pi i%
}{\tau }t}\right)  \label{ya_proc2}
\end{equation}%
Then the series has an improved convergence rate, for the Nth partial sum, of order $O\left(\sum\limits_{k=N+1}^{\infty }k^{-5}\right)
=O(N^{-4})$ \cite{kerr2022new}. This specific convergence rate assumes that the polynomials in $P_e(t)$ are of order 4. The convergence rate can be further enhanced by including additional polynomials of higher order. 

\section{New modified Laplace-Fourier method}\label{sect3}

For the solution of the NDDE (\ref{NDDE 1}) one gets that the portion of $y(t)$ related to these complex poles is given by
\begin{equation}
y(t)=2\sum_{k=1}^{\infty }\mathfrak{Re}\left( c_{k}\,e^{r_{k}t}\right),
\label{yt asy}
\end{equation}
in which, recall
\begin{equation}
c_{k}(s)=\dfrac{H(0)-bH(-\tau )+(c+bs)e^{-s\tau }\int_{-\tau
}^{0}H(v)e^{-sv}dv}{1+(bs\tau -b+c\tau )\,e^{-s\tau }}.
\label{ck asy}
\end{equation}
We also know, however, that at each of the poles $D(s)=0$. In which case 
\begin{equation}
e^{-s\tau }=\frac{s-a}{bs+c}.  \label{new imp}
\end{equation}
Thus enabling us to rewrite the denominator in Eq. (\ref{ck asy}) as
\begin{equation}
D^{\prime }(s)=1+(bs\tau -b+c\tau )\,\frac{s-a}{bs+c},
\label{Dk asy}
\end{equation}
provided $s\in \left\{ r_{k}, k\in \mathbb{N}\right\}=W$. Now applying integration by parts to the integral in the numerator of equation (\ref{ck asy}) we find that
\begin{equation}
\int_{-\tau }^{0}\!H(v)e^{-sv}dv\!=\!\frac{H(-\tau )e^{s\tau }-\!H(0)}{s}%
+\frac{H^{\prime }(-\tau )e^{s\tau }-\!H^{\prime }(0)}{s^{2}}+O\left(
s^{-3}\right).  \label{IBP}
\end{equation}
Substituting the above expression into the numerator in Eq. (\ref{ck asy}) and using Eq. (\ref{new imp}) to replace all of the exponential terms, it can be shown that
\begin{equation}
N(s)=\frac{1}{s}\left[ cH(-\tau )+aH(0)+bH^{\prime }(-\tau
)-H^{\prime }(0)\right] +O\left( s^{-2}\right),\quad s\in W.  \label{Nsk 1}
\end{equation}
Therefore 
\begin{equation}
\lim_{s\rightarrow \infty }\frac{N(s)}{D^{\prime }(s)}=\frac{1}{%
\tau s^{2}}\left[ cH(-\tau )+aH(0)+bH^{\prime }(-\tau )-H^{\prime }(0)\right] +O\left( s^{-3}\right),\quad s\in W.  \label{limit asy ck}
\end{equation}
Recall, however, that the Fourier series results that are used require an asymptotic expansion with respect to the frequencies as $k\rightarrow \infty$. Therefore, substituting $s_{k}=\ln (b)/\tau +i\alpha _{k}$ into the above
equation one gets
\begin{equation}
c_{k}^{a}=\lim_{k\rightarrow \infty }\frac{N(s_{k})}{D^{\prime }(s_{k})}=
\frac{-a_{2}}{\alpha _{k}^{2}}+i\frac{a_{3}}{\alpha _{k}^{3}}+O\left( \alpha
_{k}^{-4}\right)  \label{limit asy alpha k}
\end{equation}
where 
\begin{equation}
a_{2}=\frac{1}{\tau }\left[ cH(-\tau )+aH(0)+bH^{\prime }(-\tau )-H^{\prime
}(0)\right]  \label{a2 asy}
\end{equation}
and
\begin{align}
a_{3}=&\frac{1}{\tau }\left[ (a-2\ln (b)/\tau )(cH(-\tau )+aH(0)+bH^{\prime
}(-\tau )-H^{\prime }(0))+cH^{\prime }(-\tau )\right]\notag\\
&+\frac{1}{\tau }\left[aH^{\prime }(0)+bH^{\prime\prime }(-\tau )-H^{\prime \prime }(0)\right].  \label{a3 asy}
\end{align}

The original Laplace-Fourier method uses the formula (\ref{sk even}) for the poles $s_{k}$ and substitutes them 
into Eq. (\ref{ck asy}) right at the outset. As
a consequence all of the exponential terms simplify to
\begin{equation}
e^{-s_{k}\tau }=e^{-(\ln (b)/\tau +i\alpha _{k})\tau }=\frac{1}{b}, \qquad \forall k.
\label{exp pap2}
\end{equation}
With, observe, zero imaginary part. Therefore the term,
\begin{equation}
s_{k}\,e^{-s_{k}\tau }=\frac{1}{b}\left(\dfrac{\ln (b)+2k\pi i}{\tau }\right)
\label{Dp pap2}
\end{equation}
However the actual/exact pole $r_{k}$ is only relatively close to $s_{k}$. Therefore the term $e^{-r_{k}\tau }$ has a relatively small non-zero imaginary part. Let's denote this imaginary part as  $\epsilon_{k}i$. Then the formula for $\epsilon_{k}$ can be derived via the asymptotic expansion for $\frac{s-a}{bs+c},$ which is given by $\frac{1}{b}(1-(a+c/b)s^{-1}+O(s^{-2}))$. Thus, when multiplied by $r_{k}$ one gets
\begin{equation}
r_{k}e^{-r_{k}\tau }=\frac{1}{b}\left(\dfrac{\ln (b)+2k\pi i}{\tau}-(a+c/b)\right) +O(k^{-1}),
\end{equation}
for $k$ large. In order to distinguish the actual poles from the approximates given by Eq. (\ref{sk even}), we shall use $r_k$. Notice that this $se^{-s\tau }$ term appears in both the numerator and the denominator of Eq. (\ref{ck asy}). Therefore, the expansion of $c_{k}$ in the original Laplace-Fourier method is only exact when $ab+c=0$. However, when the contribution from the $a+c/b$ term is significant the original Laplace-Fourier method is oftentimes only marginally better than the pure Laplace solution, and in some extreme cases might be less accurate. In this paper, we have developed a new modified Laplace-Fourier method that addresses this shortcoming. In the next section we will analyze the convergence of the new modified Laplace-Fourier method.

\section{Convergence rate of the new modified Laplace-Fourier solution}

In order to establish the convergence rate of the series in Eq. (\ref{ya_proc}) we shall need an improved approximate for the location of the poles for $k$ large. This can be obtained by assuming that $s_{k}$ is of the form
\begin{equation}
s_{k}=\frac{\ln (b)+2k\pi i}{\tau }+\frac{\beta }{k^{2}}+i\frac{\gamma }{k}, \quad k\in \mathbb{N}  \label{sk imp}
\end{equation}
The constants $\beta $ and $\gamma$ can be determined by substituting Eq. (\ref{sk imp})  into Eq. (\ref{new imp}), then equating the coefficients in the respective asymptotic expansions for $k$ large. For the term on the
left it can be shown that
\begin{equation}
e^{-s_{k}\tau }=\frac{1}{b}-i\frac{\gamma \tau }{bk}-\frac{\beta \tau +\frac{%
1}{2}\gamma ^{2}\tau ^{2}}{bk^{2}}+O(k^{-3}),  \label{lhs asy}
\end{equation}
and for the term on the right hand side of Eq. (\ref{new imp}) one obtains
\begin{equation}
\frac{s_{k}-a}{bs_{k}+c}=\frac{1}{b}-i\frac{\left( a+c/b\right) \tau }{2b\pi
k}-\tau ^{2}\frac{(ab+c)(b\ln (b)/\tau +c)}{4b^{3}\pi ^{2}k^{2}}+O(k^{-3}).
\label{rhs asy}
\end{equation}
Then equating the relevant coefficients one gets the following system of
equations%
\begin{equation}
\frac{\gamma \tau }{b}=\frac{\left( a+c/b\right) \tau }{2b\pi },
\label{eq1 ab}
\end{equation}
and%
\begin{equation}
\frac{\beta \tau +\frac{1}{2}\gamma ^{2}\tau ^{2}}{b}=\tau ^{2}\frac{%
(ab+c)(b\ln (b)/\tau +c)}{4b^{3}\pi ^{2}}. \label{eq2 ab}
\end{equation}
Solving this system for $\beta $ and $\gamma$  we find that
\begin{equation}
\beta =\frac{(ab+c)(2b\ln (b)-\tau (ab-c))}{8b^{2}\pi ^{2}},  \label{alp soln}
\end{equation}
and 
\begin{equation}
\gamma =\frac{-\left( ab+c\right) }{2b\pi}. \label{bta soln}
\end{equation}
These results can now be used to verify that the convergence rate of the
new improved Laplace-Fourier solution is $O(N^{-3})$. In order to accomplish this let us recall Eq. (\ref{ya_proc}), 
\begin{equation}
y(t)\simeq c_{0}e^{r_{0}t}+2e^{\frac{\ln (b)}{\tau }t}P_{e}(t)+2%
\sum_{k=1}^{N}\mathfrak{Re}\left( c_{k}e^{r_{k}t}-c_{k}^{a}e^{\left[ \frac{\ln
(b)+2k\pi i}{\tau }\right] t}\right)  \label{ya proc 2}
\end{equation}
and focus on the residues in the sum. Now, we know that for $s$ large
\begin{equation}
c_{k}(s)=\frac{b_{2}}{s^{2}}+\frac{b_{3}}{s^{3}}+\frac{b_{4}}{s^{4}}+O\left( s^{-5}\right), \quad
\quad s\in W,  \label{asy exp s}
\end{equation}
in which
\begin{equation}
b_{2}=\frac{1}{\tau }\left[ cH(-\tau )+aH(0)+bH^{\prime }(-\tau )-H^{\prime
}(0)\right],  \label{b2 asy}
\end{equation}
and
\begin{equation}
b_{3}=\frac{1}{\tau }\left[ a^{2}H(0)+acH(-\tau )+aH(0)+(ab+c)H^{\prime
}(-\tau )+bH^{\prime \prime }(-\tau )-H^{\prime \prime }(0)\right].
\label{b3 asy}
\end{equation}
If we use the improved approximate $s_k$ in Eq. (\ref{sk imp}) as a surrogate for $r_{k}$, then one gets
\begin{equation}
c_{k}(r_{k})=\frac{-b_{2}}{\alpha _{k}^{2}}+i\frac{b_{3}-2b_{2}\ln (b)/\tau 
}{\alpha _{k}^{3}}+\frac{3b_{2}(\ln (b)/\tau )^{2}-3b_{3}\ln (b)/\tau
+b_{4}-2b_{2}(ab+c)/\tau }{\alpha _{k}^{4}}+O\left( k^{-5}\right).
\label{ck rk}
\end{equation}
However, recall that for the second residue in our sum we must use the $s_{k}$ given by Eq. (\ref{sk even})$.$ In which case
\begin{equation}
c_{k}^{a}(s_{k})=\frac{-b_{2}}{\alpha _{k}^{2}}+i\frac{b_{3}-2b_{2}\ln
(b)/\tau }{\alpha _{k}^{3}}+\frac{3b_{2}(\ln (b)/\tau )^{2}-3b_{3}\ln
(b)/\tau +b_{4}}{\alpha _{k}^{4}}+O\left( k^{-5}\right)   \label{cka sk}
\end{equation}
In which case, for $k$ large 
\begin{equation}
c_{k}e^{r_{k}t}-c_{k}^{a}e^{\left[ \frac{\ln (b)+2k\pi i}{\tau }\right]
t}\simeq \left( c_{k}-c_{k}^{a}\right) e^{r_{k}t}=\left( -\frac{%
2b_{2}\,(ab+c)/\tau }{\alpha _{k}^{4}}+O\left( k^{-5}\right) \right) e^{r_{k}t}
\label{k - 4}
\end{equation}
As a consequence, when $ab+c\neq 0$ the convergence rate of the new modified Laplace-Fourier solution has a remarkable order of convergence $O(N^{-3})$. In the next section we will present some examples that numerically illustrate the theoretical results obtained in this work.

\section{Illustrative results of the new modified Laplace-Fourier method} \label{section4}

In this section we will present a comparison between the the pure Laplace and the Laplace-Fourier methods by solving several NDDEs. We shall show that the developed method improves the accuracy of the Laplace solution without increasing the number of terms in the truncated series. In most of the examples we will use the approximate formula for the location of the poles as our initial guess for determining the actual poles. In all of our examples, we include a graph of the analytic solution, obtained via the MoS. In order to have fair comparisons between the methods we use the  same number of terms in each of the truncated series related to the methods. In addition, we include the error graphs for the solutions generated by the pure Laplace solution and the Laplace-Fourier methods. The errors are obtained by computing the difference between the approximated solutions and the analytical solution, over some finite time interval. 
We use Maple to perform all the required computations for the Laplace transform method, Laplace-Fourier method, and the MoS. We compute the coefficients in the asymptotic expansions using Maple's series command. 

\subsection{Example 1}  
\begin{figure}[ht] 
\centering
\includegraphics[scale=0.75]{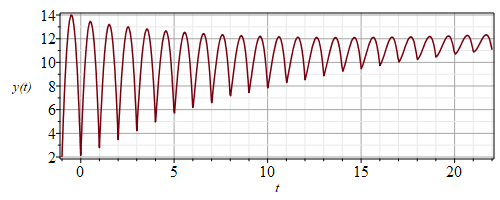}
\caption{Analytical solution obtained by the MoS. The NDDE is $y^{\prime }(t)\!=\!-2.1y(t)\!+\!0.9y^{\prime }(t-1)\!+\!2.12y(t-1), H(t)\!=\!2-48t(1+t),
t\in \lbrack -1,0]$.}
\label{fig1a}		
\end{figure}

\begin{figure}[ht] 
\centering
\includegraphics[scale=0.25]{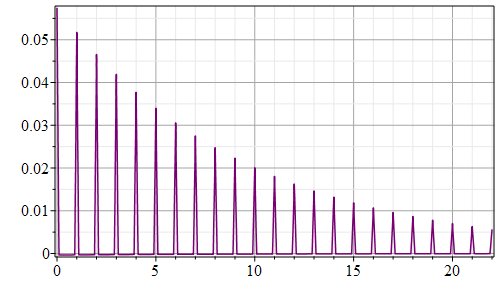}\includegraphics[scale=0.25]{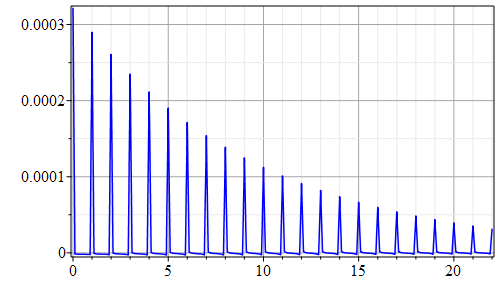}\includegraphics[scale=0.25]{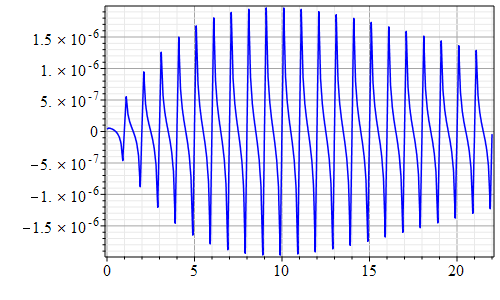}
\caption{Errors of the solutions obtained by pure Laplace (left), original Laplace-Fourier (center) and the new modified Laplace-Fourier (right) methods. The NDDE is $y^{\prime }(t)\!=\!-2.1y(t)\!+\!0.9y^{\prime }(t-1)\!+\!2.12y(t-1), H(t)\!=\!2-48t(1+t),
t\in \lbrack -1,0]$.}
\label{fig1b}		
\end{figure}

\begin{figure}[ht] 
\centering
\includegraphics[scale=0.4]{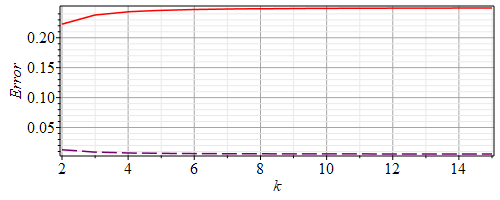}\includegraphics[scale=0.4]{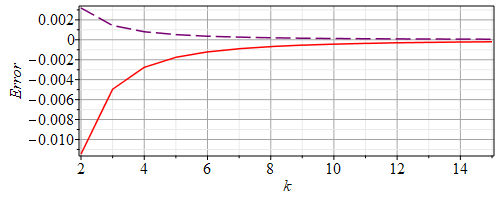}
\caption{Relative errors of the imaginary (solid lines) and real (dashed lines) parts using the correct residues $c_{k}$ versus the approximate values $c_{k}^{a}$ for the original Laplace-Fourier (left) and for the new modified Laplace-Fourier method (right). The NDDE is $y^{\prime }(t)\!=\!-2.1y(t)\!+\!0.9y^{\prime }(t-1)\!+\!2.12y(t-1), H(t)\!=\!2-48t(1+t),
t\in \lbrack -1,0]$.}
\label{fig1c}		
\end{figure}

Let us consider the following linear NDDE, 
\begin{equation}
y^{\prime }(t)\!=\!-2.1y(t)\!+\!0.9y^{\prime }(t-1)\!+\!2.12y(t-1), H(t)\!=\!2-48t(1+t),
t\in \lbrack -1,0].\!\label{ex1}
\end{equation}
We computed the solutions of NDDE (\ref{ex1}) by the MoS (piecewise exact analytical solution), pure Laplace, original Laplace-Fourier and the new modified Laplace-Fourier methods using the symbolic software Maple. For NDDE (\ref{ex1}) there are two real poles, one negative ($r_0$) and one positive ($r_1$). Fig. \ref{fig1a} shows the (oscillatory) solution obtained via the MoS. The solution may appear to approach a steady state, but it doesn't. It approaches to $c_1\,e^{r_1\,t}$, where $r_1\approx0.009$ and $c_1\approx 9.64$. Therefore, the solution of NDDE (\ref{ex1}) increases exponentially and there is no steady state. In Fig. \ref{fig1b}, the graphs of the absolute errors related to the different methods are presented. It is clear that the solution given by the new modified Laplace-Fourier method is the most accurate. While, as expected, the solution produced by the standard pure Laplace method is the least accurate. The accuracy can be further enhanced by increasing the number of terms in the series and/or by including more (higher degree) polynomials. The errors for all solutions decrease as time increases, despite the fact that the solution itself is increasing. This is an important feature of all the implemented methods. For the absolute error plots, the peaks which occur at $t=k\tau$ are due to the nature of the non-harmonic series form of all the solutions, generated by these methods \cite{kerr2022accuracy,kerr2022new}. It is important to remark, that for all these methods once the solution is obtained we can compute the solution for any time, even for a very large value of $t$, with a single calculation. In this example polynomials of degree 8 were used. We note, however, that for this particular example, the number of polynomials has a minimal effect in reducing the errors once we go beyond $p_3(x)$. 
For the NNDE (\ref{ex1}), we have that $ab+c=0.23$. On the left hand side of Fig. \ref{fig1c} we have plotted the relative errors of the real and imaginary parts computed using the correct residues $c_{k}$ versus the approximate values $c_{k}^{a}$ which are obtained from using the original Laplace-Fourier method. In particular, the graphs of $\mathfrak{Re}(c_{k}-c_{k}^{a})/\mathfrak{Re}(c_{k})$ and $\mathfrak{Im}(c_{k}-c_{k}^{a})/
\mathfrak{Im}(c_{k})$ for $k=2 \ldots 12$. On the right hand side of Fig. \ref{fig1c} we have plotted the
analogous errors which are obtained from using the new modified Laplace-Fourier approximate values $c_{k}^a$ given by Eq. (\ref{limit asy alpha k}). 

\subsection{Example 2}  
\begin{figure}[ht] 
\centering
\includegraphics[scale=0.75]{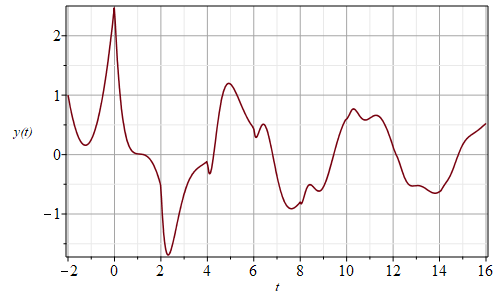}
\caption{Analytical solution obtained by the MoS. The NDDE is $y^{\prime }(t)\!\!=\!\!-2.1y(t)\!+\!\frac{7}{11}y^{\prime }(t-2)\!-\!2y(t-2),H(t)\!=\!1\!+\!\frac{3}{2}(t+2)(0.5+t), t\in [-2,0]$.}
\label{fig2a}		
\end{figure}

\begin{figure}[ht] 
\centering
\includegraphics[scale=0.27]{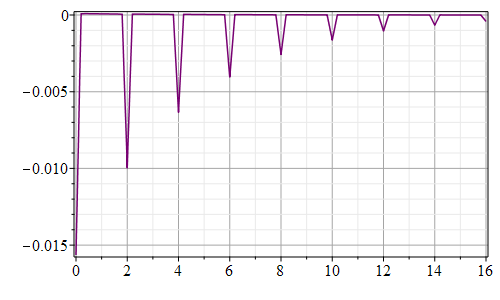}\includegraphics[scale=0.27]{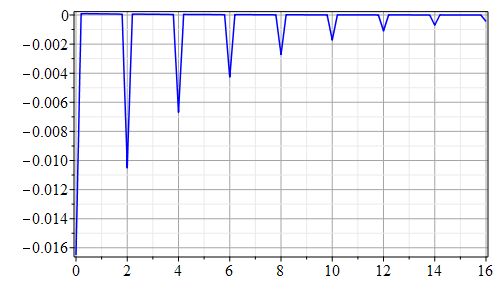}\includegraphics[scale=0.27]{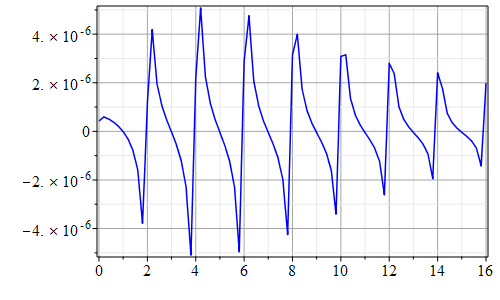}
\caption{Errors of the solutions obtained by pure Laplace (left), original Laplace-Fourier (center) and the new modified Laplace-Fourier (right) methods. The NDDE is $y^{\prime }(t)\!\!=\!\!-2.1y(t)\!+\!\frac{7}{11}y^{\prime }(t-2)\!-\!2y(t-2),H(t)\!=\!1\!+\!\frac{3}{2}(t+2)(0.5+t), t\in [-2,0]$.}
\label{fig2b}		
\end{figure}

Let us consider the following linear NDDE, 
\begin{equation}
y^{\prime }(t)\!\!=\!\!-2.1y(t)\!+\!\frac{7}{11}y^{\prime }(t-2)\!-\!2y(t-2),H(t)\!=\!1\!+\!\frac{3}{2}(t+2)(0.5+t), t\in [-2,0].\label{ex2}
\end{equation}
We computed the solutions of NDDE (\ref{ex2}) by the MoS (piecewise exact analytical solution), pure Laplace, original Laplace-Fourier and the new modified Laplace-Fourier methods. For NDDE (\ref{ex2}) there are no real poles. In Fig. \ref{fig2a} 
we show the solution computed by the MoS. In Fig. \ref{fig2b}, the graphs of the absolute errors related to the different methods are presented. Again, it is clear that the solution given by the new modified Laplace-Fourier method is the most accurate. In this example, however, the solution produced by the original Laplace-Fourier method is slightly less accurate. This is due to the fact that there is a significant error in the asymptotic expansion for the residues $c_k$. However, the solution generated by the new modified Laplace-Fourier method is highly accurate due to the improvement in the computation of the residues, which is one of the main
results of this paper. As usual, more accuracy can be obtained if the number of terms in the solutions are increased and/or more polynomials are included. The errors for all solutions decrease as time increases \cite{kerr2022new,kerr2022accuracy}. In this example polynomials of degree 8 were used. 

\begin{figure}[ht] 
\centering
\includegraphics[scale=0.4]{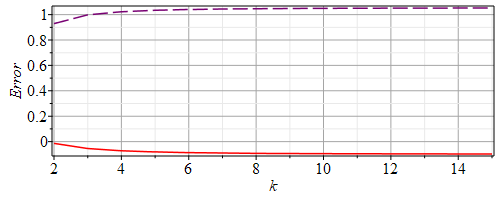}\includegraphics[scale=0.4]{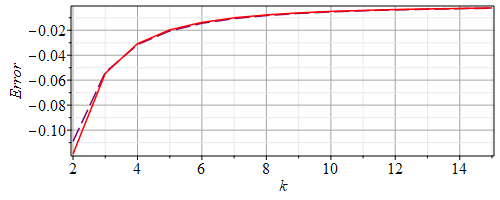}
\caption{Relative errors of the imaginary (solid lines) and real (dashed lines) parts using the correct residues $c_{k}$ versus the approximate values $c_{k}^{a}$ for the original Laplace-Fourier (left) and for the new modified Laplace-Fourier method (right). The NDDE is $y^{\prime }(t)\!=\!-2.1y(t)\!+\!0.9y^{\prime }(t-1)\!+\!2.12y(t-1), H(t)\!=\!2-48t(1+t),
t\in \lbrack -1,0]$.}
\label{fig2c}		
\end{figure}

For NNDE (\ref{ex2}), we have that $ab+c=-367/110$. On the left hand side of Fig. \ref{fig2c} we have plotted the relative errors of the real and imaginary parts computed using the correct residues $c_{k}$ versus the approximate values $c_{k}^{a}$ which are obtained from using the original Laplace-Fourier method. In particular, the graphs of $\mathfrak{Re}(c_{k}-c_{k}^{a})/\mathfrak{Re}(c_{k})$ and $\mathfrak{Im}(c_{k}-c_{k}^{a})/
\mathfrak{Im}(c_{k})$ for $k=2 \ldots 12$. On the right hand side of Fig. \ref{fig2c} we have plotted the
analogous errors which are obtained from using the new modified Laplace-Fourier approximate values $c_{k}^a$ given by Eq. (\ref{limit asy alpha k}).

\subsection{Example 3}  
\begin{figure}[ht] 
\centering
\includegraphics[scale=0.5]{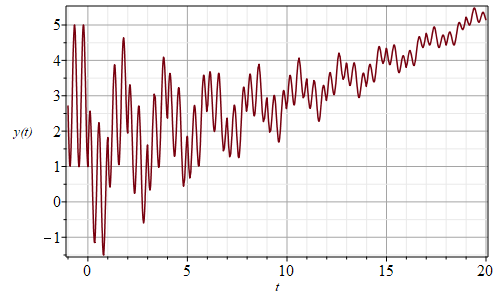}
\caption{Analytical solution obtained by the MoS. The NDDE is $y^{\prime }(t)\!\!=\!\!\frac{14}{33}y(t)\!-\!\frac{8}{9}y^{\prime }(t-1)\!-\!\frac{1}{3}y(t-1),H(t)\!=\!3-2\cos 14t, t\in [-1,0]$.}
\label{fig3a}		
\end{figure}

\begin{figure}[ht] 
\centering
\includegraphics[scale=0.3]{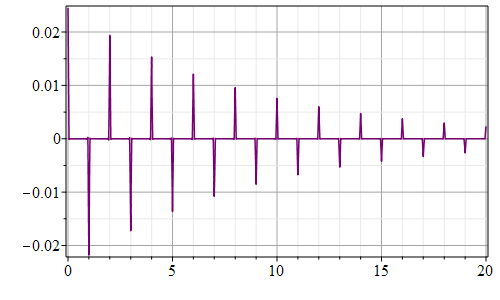}\includegraphics[scale=0.3]{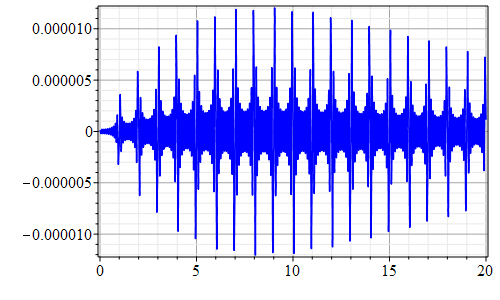}
\caption{Errors of the solutions obtained by pure Laplace (left) and the new modified Laplace-Fourier (right) methods. The NDDE is $y^{\prime }(t)\!\!=\!\!\frac{14}{33}y(t)\!-\!\frac{8}{9}y^{\prime }(t-1)\!-\!\frac{1}{3}y(t-1),H(t)\!=\!3-2\cos 14t, t\in [-1,0]$. The number of terms used to compute the approximated solution (\ref{ya_proc2}) is $N=50$.}
\label{fig3b}		
\end{figure}

The last NDDE that we consider includes a periodic history function and also considers a coefficient $b<0$. Recall that for this particular case the formulas for the approximate poles, frequencies, polynomials etc. are different (\cite{kerr2022new}). The NDDE is the following 
\begin{equation}
y^{\prime }(t)\!\!=\!\!\frac{14}{33}y(t)\!-\!\frac{8}{9}y^{\prime }(t-1)\!-\!\frac{1}{3}y(t-1),H(t)\!=\!3-2\cos 14t, t\in [-1,0].\label{ex3}
\end{equation}
We compute the analytical solution of NDDE (\ref{ex3}) in order to compare it with the solutions provided by the pure Laplace and the new modified Laplace-Fourier methods. For NDDE (\ref{ex3}) there is one real pole (positive). Fig. \ref{fig3a} shows the analytical solution given by the MoS. In Fig. \ref{fig3b}, the graphs of the absolute errors related to the pure Laplace and the new modified Laplace-Fourier methods are presented. It can be seen that the solution given by the new modified Laplace-Fourier method is significantly more accurate. This result is expected due to the fact that the new modified Laplace-Fourier method incorporates the improved asymptotic expansion for the residues.

\begin{figure}[ht] 
\centering
\includegraphics[scale=0.3]{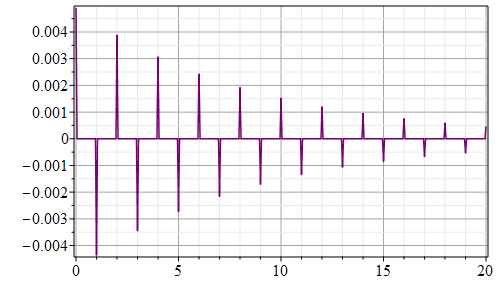}\includegraphics[scale=0.3]{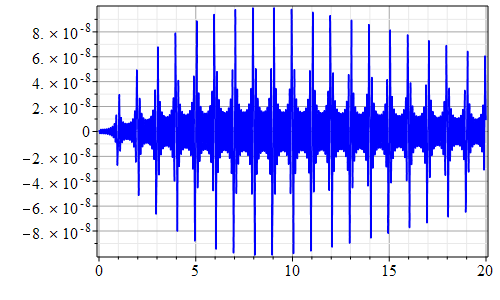}
\caption{Relative errors of the imaginary (left) and real (right) parts using the correct residues $c_{k}$ versus the approximate values $c_{k}^{a}$ for the original Laplace-Fourier (solid lines) and for the new modified Laplace-Fourier method (dashed lines). The NDDE is $y^{\prime }(t)\!\!=\!\!\frac{14}{33}y(t)\!-\!\frac{8}{9}y^{\prime }(t-1)\!-\!\frac{1}{3}y(t-1),H(t)\!=\!3-2\cos 14t, t\in [-1,0]$. The number of terms used to compute the approximated solution (\ref{ya_proc2}) is $N=250$.}
\label{fig3c}		
\end{figure}

In this example we also show an important feature of the new modified Laplace-Fourier. In Fig. \ref{fig3c}, we have plotted the error graphs for the pure Laplace and the new modified Laplace-Fourier solutions, but in
this case we increased the number of terms for the approximated solution (\ref{ya_proc2}) to $N=250$. Notice, that the improvement in the new modified Laplace-Fourier method is considerably larger than the one of the pure Laplace method. We also note that Figs. \ref{fig3b} and \ref{fig3c} provide corroboration that the convergence rate for the new improved Laplace-Fourier method is $O(N^{-3})$. We know the Pure Laplace convergence rate is $O(N^{-1})$. And we increased the number of terms in the series by a factor of 5. If we measure the corresponding error reduction based on the y-axis scales, we have an error reduction factor for pure Laplace $\approx 0.02/0.004=5$. However, the analogous error reduction factor for the improved Fourier-Laplace $\approx0.000012/(9.6\times10^{-8}) = 125$. In Table \ref{tab1} we present the maximum errors of the pure Laplace and the modified Laplace-Fourier methods. We have also, however, included the error for $N = 500$. It is clear that the modified Laplace-Fourier method has a much better convergence rate than the pure Laplace method. In this example polynomials of degree 7 were used. And, although not included here, similar results were also observed for the NDDEs (\ref{ex1}) and (\ref{ex2}). Notice that once the modified Laplace-Fourier solution is obtained it can evaluated at any value of $t$ with a very low computational time, which is a main advantage in comparison with numerical methods.


\begin{center}
\begin{tabular}{|p{1cm}|p{2.5cm}|p{2.2cm}|p{3cm}|p{2.2cm}|}
\hline
\multicolumn{5}{|c|}{Error} \\
\hline
N & Pure Laplace &Convergence & Laplace-Fourier & Convergence \\
\hline
50  & $\approx 0.02$ & & $\approx 1.2\times10^{-5}$ & \\
\hline
250 & $\approx 0.004$ & $O(N^{-1})$   & $\approx 9.6\times10^{-8}$ & $O(N^{-3})$ \\
\hline 
500 & $\approx 0.002$ & $O(N^{-1})$&  $\approx 1.2\times10^{-8}$ & $O(N^{-3})$ \\
\hline
\end{tabular}\label{tab1}
\end{center}

\section{Conclusions}\label{section6}

In this article, we have developed a new modified Laplace-Fourier method to obtain solutions of linear NDDEs. This new modified method still combines the Laplace transform and Fourier series theory to obtain the solutions. However, by deriving and implementing an improved formula for the asymptotic expansion, for the residues, we were able to obtain more accurate solutions. The newly developed modified Laplace-Fourier method generates more accurate solutions than the ones generated by the pure Laplace method and the original Laplace-Fourier method. We have shown that the convergence rate of the new modified Laplace-Fourier solution has a remarkable order of convergence $O(N^{-3})$. Examples and different comparisons with other methods have illustrated the accuracy of the new developed method. We have implemented the new modified method using the well-known Maple computer algebra software, but it can also be implemented with other available symbolic software such as Mathematica or Maxima. The solutions generated by the new modified Laplace-Fourier method were very close to the analytical ones, with impressively small errors. Moreover, the errors effectively become negligible as $t$ increases. It is important to remark that all of the advantageous features from the original Laplace-Fourier method are retained. For instance, the solution still accounts for the terms in the tail of the series, but now in a much more effective manner. And the accuracy of the solution can be increased by increasing the number of terms in the series and/or increasing the degree of the polynomials in the piecewise function.
Furthermore, the modified Laplace-Fourier method generates a solution over the entire time domain, thus enabling one to compute the solution at any time with a
single calculation. This is a main advantage in comparison with numerical methods. Finally, the solution given by the new modified Laplace-Fourier method improves at a much larger rate than the one of the pure Laplace method. This, is an excellent additional advantage of the new modified Laplace-Fourier. The validity of the modified technique is corroborated by means of various illustrative examples. In summary, the modified Laplace-Fourier method is a highly efficient method for solving linear NDDEs. 

%
\section*{Conflict of interest}

The authors declare that they have no conflict of interest.




%
%

\end{document}